\documentclass[12pt,reqno]{amsart}

\newcommand\version{September 6, 2021}

\usepackage{amsmath,amsfonts,amsthm,amssymb,amsxtra}
\usepackage{bbm} 

\newtheorem{theorem}{Theorem}

\newtheorem{problem}[theorem]{Problem}

\theoremstyle{definition}

\theoremstyle{remark}

\newcommand{\1}{\mathbbm{1}}

\renewcommand{\epsilon}{\varepsilon}

\newcommand{\f}{\mathrm{f}}
\newcommand{\gld}{\mathrm{gld}}
\newcommand{\gKS}{\mathrm{gKS}}

\newcommand{\N}{\mathbb{N}}

\renewcommand{\phi}{\varphi}
\newcommand{\R}{\mathbb{R}}
\newcommand{\rgKS}{\mathrm{rgKS}}

\newcommand{\Sph}{\mathbb{S}}

\DeclareMathOperator{\per}{Per}

\begin{document}

\title[Some minimization problems --- \version]{Some minimization problems\\ for mean field models with competing forces}

\author{Rupert L. Frank}
\address[Rupert L. Frank]{Mathe\-matisches Institut, Ludwig-Maximilans Universit\"at M\"unchen, The\-resienstr.~39, 80333 M\"unchen, Germany, and Munich Center for Quantum Science and Technology, Schel\-ling\-str.~4, 80799 M\"unchen, Germany, and Mathematics 253-37, Caltech, Pasa\-de\-na, CA 91125, USA}
\email{r.frank@lmu.de}

\renewcommand{\thefootnote}{${}$} \footnotetext{\copyright\, 2021 by the author. This paper may be reproduced, in its entirety, for non-commercial purposes.}

\begin{abstract}
	We review recent results on three families of minimization problems, defined on subsets of nonnegative functions with fixed integral. The competition between attractive and repulsive forces leads to transitions between parameter regimes, where minimizers exist and where they do not. The problems considered are generalized liquid drop models, swarming models and generalized Keller--Segel models.
\end{abstract}

\maketitle

\section{Introduction}

In this survey we discuss three families of minimization problems. They are simple mathematical toy models for physical or biological phenomena. While their origins are rather different, they share some mathematical similarities and differences and we think it is worthwhile to look at them side by side.

The common feature of all three problems is that they are of mean-field type. They involve an `energy' functional that is defined on a subset of nonnegative functions (`densities') whose integral is fixed (`total mass'). They are, at least on a heuristic level, derived from microscopic, many-body models. The densities in the mean-field models describe the distribution of the microscopic particles in the limit of a large number of particles, and similarly the energy functionals in our models are obtained as macroscopic approximations to microscopic energy functionals.

Another common feature of the problems discussed here is that the energy functionals have two contributions that compete with each other. There are attractive forces that keep the particles together and try to concentrate them and there are repulsive forces that push them apart and try to spread them out. Typically, these forces act on different length scales and one is of short range and the other one of long range type. The existence of a minimizer can be understood as the forces being in a local equilibrium, while the nonexistence typically means that one of the forces dominates the other.

We are particularly interested in situations where, as a parameter of the problem is varied continuously, there is either a transition between existence and nonexistence of minimizers, or a sharp change in the properties of minimizers. A typical parameter that is varied is the total mass, but in one of the models it is also a parameter describing the shape of the forces acting between the particles.

\subsection*{The models}

Let us be more specific about the three families of models that we will consider. Throughout, $N\geq 1$ is the dimension of the underlying Euclidean space.

For the \emph{generalized liquid drop model}, depending on a parameter $\lambda\in (0,N)$, we define for any measurable set $\Omega\subset\R^N$,
\begin{equation}
	\label{eq:gld}
	\mathcal E_\lambda^\gld[\Omega] := \per \Omega + \frac12 \iint_{\Omega\times\Omega} \frac{dx\,dy}{|x-y|^\lambda} \,.
\end{equation}
Here $\per \Omega$ denotes the perimeter in the sense of De Giorgi; see, e.g., \cite{Ma}. The corres\-ponding minimization problem is, for $m\in(0,\infty)$,
\begin{equation}
	\label{eq:gldmin}
	E_\lambda^\gld(m) := \inf\left\{ \mathcal E_\lambda^\gld[\Omega]:\ \Omega\subset\R^N \ \text{measurable},\ |\Omega|= m \right\}.
\end{equation}
The original liquid drop model, suggested by Gamow \cite{Ga} for the description of atomic nuclei, corresponds to $\lambda=1$ in dimension $N=3$.

For the \emph{flocking model}, depending on parameters $\lambda\in(0,N)$ and $\alpha\in (0,\infty)$, we define for any nonnegative, measurable function $\rho$ on $\R^N$,
\begin{equation}
	\label{eq:f}
	\mathcal E_{\lambda,\alpha}^\f[\rho] := \frac12 \iint_{\R^N\times\R^N} \rho(x) \left( |x-y|^{-\lambda} + |x-y|^\alpha \right)\rho(y)\,dx\,dy \,. 
\end{equation}
The corresponding minimization problem is, for $m\in(0,\infty)$,
\begin{equation}
	\label{eq:fmin}
	E_{\lambda,\alpha}^\f(m) := \inf\left\{ \mathcal E_{\lambda,\alpha}^\f[\rho]:\ \rho\in L^1(\R^N),\ 0 \leq\rho\leq 1 \,,\ \int_{\R^N} \rho\,dx = m \right\}.
\end{equation}
This model was suggested by Burchard, Choksi and Topaloglu \cite{BCT}. It is a simple model to describe the flocking behavior in stable states of a large group of animals such as fish or birds.

For the \emph{generalized Keller--Segel model}, depending on parameters $q\in(0,1)$ and $\alpha\in (0,\infty)$, we define for any nonnegative function $\rho\in L^q(\R^N)$,
\begin{equation}
	\label{eq:gks}
	\mathcal E_{q,\alpha}^\gKS[\rho] := - \int_{\R^N} \rho^q\,dx + \frac12 \iint_{\R^N\times\R^N} \rho(x) |x-y|^\alpha \rho(y)\,dx\,dy \,.
\end{equation}
The corresponding minimization problem is
\begin{equation}
	\label{eq:gksmin}
	E_{q,\alpha}^\gKS := \inf\left\{ \mathcal E_{q,\alpha}^\gKS[\rho]:\ 0 \leq \rho\in L^q(\R^N),\ \int_{\R^N} \rho\,dx = 1 \right\}.
\end{equation}
Note that here, in contrast to the two previous problems, we fix the integral of $\rho$ to be one. The more general case, where it is fixed to be equal to $m$, can be reduced to the present one by scaling. The generalized Keller--Segel model was introduced in \cite{CDDFH} and generalizes the standard Keller--Segel model, which corresponds (after some rescaling) to the limit cases $q=1$ and $\alpha=0$ in dimension $N=2$.

\subsection*{Competing forces}

Let us discuss in which sense in the above models two forces compete with each other.

In the generalized liquid drop model, the perimeter term corresponds to an attractive short range force, whereas the double integral term corresponds to a repulsive long range force. Note that by the isoperimetric inequality, see, e.g., \cite{Ma},
$$
\inf\left\{ \per \Omega:\ \Omega\subset\R^N \ \text{measurable},\ |\Omega|=m \right\} = N^\frac{N-1}{N} |\Sph^{N-1}|^\frac{1}{N} m^\frac{N-1}{N}
$$
with equality if and only if $\Omega$ is a ball (up to sets of measure zero). On the other hand, it is easy to see that
$$
\inf\left\{ \frac{1}{2} \iint_{\Omega\times\Omega} \frac{dx\,dy}{|x-y|^\lambda}:\ \Omega\subset\R^N \ \text{measurable},\ |\Omega|=m \right\} = 0
$$
and the infimum is not attained. A minimizing sequence is given, for instance, by taking $\Omega$ as a union of a large number of small balls placed very far apart from each other. Next, we note that, by scaling,
$$
E_\lambda^\gld(m) = \inf\left\{ m^\frac{N-1}{N} \per\omega + m^\frac{2N-\lambda}{N} \frac{1}{2} \iint_{\omega\times\omega} \frac{dx\,dy}{|x-y|^\lambda} :\ \omega\subset\R^N \ \text{meas.},\ |\omega|= 1 \right\}.
$$
Since $(N-1)/N<(2N-\lambda)/N$, the perimeter term is dominant for small $m$, whereas the double integral is dominant for large $m$. We therefore expect existence of minimizers for small $m$, whereas for large $m$ we might have nonexistence of minimizers.

In the flocking model, the $\alpha$-term corresponds to an attractive force, while the $\lambda$-term corresponds to a repulsive force. Moreover, the $\alpha$-term is relevant on large distances and the $\lambda$-term on short ones. By rearrangement inequalities and the bathtub principle,
$$
\inf\left\{ \frac12 \iint_{\R^N\times\R^N} \rho(x)\,|x-y|^\alpha \rho(y)\,dx\,dy:\ \rho\in L^1(\R^N),\, 0\leq \rho\leq 1,\, \int_{\R^N} \rho\,dx = m \right\}
$$
is attained if and only if $\rho$ is the characteristic function of a ball of volume $m$. Moreover, as a consequence of what we said in the generalized liquid drop model,
$$
\inf\left\{ \frac12 \iint_{\R^N\times\R^N} \frac{\rho(x)\,\rho(y)}{|x-y|^\lambda}\,dx\,dy:\ \rho\in L^1(\R^N),\, 0\leq \rho\leq 1,\, \int_{\R^N} \rho\,dx = m \right\} = 0
$$
and the infimum is not attained. Next, we note that, by scaling,
\begin{align*}
	E_{\lambda,\alpha}^\f(m) & = \inf\left\{ m^\frac{2N-\lambda}{N} \frac12 \iint_{\R^N\times\R^N} \frac{\sigma(x)\,\sigma(y)}{|x-y|^\lambda}\,dx\,dy \right. \\
	& \qquad\qquad + m^\frac{2N+\alpha}{N} \frac12 \iint_{\R^N\times\R^N} \sigma(x)\,|x-y|^\alpha \sigma(y)\,dx\,dy : \\
	& \qquad\qquad \left. \sigma\in L^1(\R^N),\, 0\leq \sigma\leq 1,\, \int_{\R^N} \sigma\,dx = 1 \right\}.
\end{align*}
Since $(2N-\lambda)/N<(2N+\alpha)/N$, the $\alpha$-term is dominant for large $m$ and we expect existence of minimizers and closeness to the characteristic function of a ball.  We also have
\begin{align*}
	E_{\lambda,\alpha}^\f(m) & = m^2 \inf\left\{ \frac12 \iint_{\R^N\times\R^N} \sigma(x)\left( |x-y|^{-\lambda} + \right)|x-y|^\alpha \right) \sigma(y)\,dx\,dy : \\
	& \qquad\qquad\qquad \left. \sigma\in L^1(\R^N),\, 0\leq \sigma\leq m^{-1},\, \int_{\R^N} \sigma\,dx = 1 \right\}.
\end{align*}
For small $m$, we expect that the constraint $\sigma\leq m^{-1}$ is irrelevant and that the minimizer is $m$ times the minimizer of the problem
\begin{align*}
	\inf\left\{ \frac12\! \iint_{\R^N\times\R^N} \!\!\!\sigma(x) \! \left( |x-y|^{-\lambda} + |x-y|^\alpha \right)\!\sigma(y)\,dx\,dy: 0\leq \sigma\in L^1(\R^N), \int_{\R^N} \!\sigma\,dx = 1 \right\}\!,
\end{align*}
provided a minimizer for the latter problem exists and is bounded.

Finally, in the generalized Keller--Segel model, the $L^q$ term corresponds to a repulsive short range force, whereas the double integral term corresponds to an attractive long range force. Note that
$$
\inf\left\{ -\int_{\R^N} \rho^q\,dx:\ 0 \leq \rho\in L^q(\R^N),\, \int_{\R^N} \rho\,dx =1 \right\} = - \infty \,.
$$
A minimizing sequence is given, for instance, by a sequence that spreads out like $\ell^{-N} \sigma(x/\ell)$ with $\ell\to\infty$. On the other hand,
$$
\inf\left\{ \frac12 \iint_{\R^N\times\R^N} \rho(x)\,|x-y|^\alpha \rho(y)\,dx\,dy:\ 0\leq \rho\in L^q(\R^N),\, \int_{\R^N} \rho\,dx = 1 \right\} = 0
$$
and the infimum is not attained. A minimizing sequence is given, for instance, by a delta sequence $\ell^{-N} \sigma(x/\ell)$ with $\ell\to 0$. Since, as we already mentioned, in this model the dependence on the total mass is trivial, we are looking here for a transition in terms of the parameters $q$ and $\alpha$. Intuitively, the repulsive force is stronger the smaller $q$ and the attractive force is stronger the larger $\alpha$. The above examples suggest that two mechanisms for the nonexistence of a minimizer are conceivable, namely both spreading out and concentration of minimizing sequences.

\subsection*{Structure of the paper}

In the following three sections we summarize what is known about the three families of minimization problems. The presentation will be rather compact and we refer to the original papers for the proofs. We do, however, emphasize several open questions concerning each model. In a short appendix we provide details for a simple, unpublished results in the one-dimensional generalized liquid drop model.


\subsection*{Acknowledgement}

The author would like to thank the organizers of the 8th European Congress of Mathematics for the organization of the meeting and for the invitation to speak. Since the topic of his invited talk was recently and rather exhaustively reviewed in \cite{FLT}, this contribution is based on a talk in a minisymposium at the congress, organized by L.~Pick, to whom the author is very grateful. The results reviewed here were obtained in collaboration with many researchers and it is a pleasure to thank, in particular, Jos\'e Carrillo, Mat\'ias Delgadino, Jean Dolbeault, Franca Hoffmann, Rowan Killip, Mathieu Lewin, Elliott Lieb and Phan Th\`anh Nam for many stimulating discussions. Partial support through U.S. National Science Foundation grants DMS-1363432 and DMS-1954995 and through the Deutsche Forschungsgemeinschaft (DFG, German Research Foundation) through Germany’s Excellence Strategy EXC - 2111 - 390814868 is acknowledged.


\section{The generalized liquid drop model}

In this section we consider the energy functional \eqref{eq:gld} and the corresponding minimization problem \eqref{eq:gldmin}. We assume throughout that $0<\lambda<N$.

\medskip

Let us set, for fixed $\lambda$ and $N$,
$$
m_* := \left( \frac{2^{1/N}-1}{1-2^{-(N-\lambda)/N}}\, \frac{\per B_1}{\frac12\iint_{B_1\times B_1}|x-y|^{-\lambda}\,dx\,dy} \right)^{N/(N-\lambda+1)} |B_1| \,,
$$
where $B_1$ denotes the unit ball in $\R^N$. The number $m_*$ is the unique solution $m>0$ of the equation
\begin{equation}
	\label{eq:charmstar}
	\mathcal E^\gld_{\lambda}\left[ \left( \frac{m}{|B_1|}\right)^{1/N} B_1 \right] = 2\, \mathcal E^\gld_{\lambda}\left[ \left( \frac{m}{2|B_1|}\right)^{1/N} B_1 \right].
\end{equation}
Thus, the energy of a ball of mass $m_*$ is equal to the energy of two balls, each of mass $m_*/2$, placed infinitely far apart. For $m<m_*$ one has $<$ instead of $=$ in \eqref{eq:charmstar} and for $m>m_*$ one has $>$.

In the physics literature it is typically taken for granted that in the special case $\lambda=1$ and $N=3$, balls are minimizers for $E^\gld_\lambda(m)$ for $m\leq m_*$ and there is no minimizer for $m>m_*$. In the mathematics literature this appears explicitly as a conjecture in work of Choksi and Peletier \cite{ChoPel-10,ChoPel-11}.

One may wonder whether the analogous conjecture is valid in the general case $0<\lambda<N$.  In dimension $N=1$ this is indeed the case, as can be verified by elementary computations; see Appendix \ref{sec:app}. It is shown in \cite{KnuMur-13,BonCri-14} that for any $N\geq 2$ there is a $\lambda_c>0$ such that for all $0<\lambda<\lambda_c$ the conjecture is true; see \cite{MurZal-14} for an explicit lower bound on $\lambda_c$ for $N=2$. In the remaining cases, the validity or invalidity of the conjecture is open.

\medskip

\emph{Existence.} As a first step towards this conjecture, before asking whether minimizers for $E^\gld_{\lambda}(m)$ are balls for all $m\leq m_*$, it is natural to ask whether minimizers exist for all $m\leq m_*$. This is indeed the case, as shown in \cite{FrNa}. Moreover, it is shown there as well that if there are no minimizers for $m>m_*$, then balls are minimizers for $m\leq m_*$.

The proof of \cite{FrNa} proceeds by verifying that for any $m<m_*$ one has the strict binding inequality
$$
E^\gld_{\lambda}(m) < E^\gld_{\lambda}(m') + E^\gld_{\lambda}(m-m')
\qquad\text{for all}\ 0<m'<m \,.
$$
According to a compactness result in \cite{FraLie-15} this implies the existence of a minimizer for $E^\gld_{\lambda}(m)$ for $m\leq m_*$.

\medskip

\emph{Uniqueness.} We adress the question of whether balls are minimizers. A convexity argument due to Bonacini and Cristoferi \cite[Theorem 2.10]{BonCri-14} shows that there is a number $m_c^{\rm ball}\in [0,\infty)\cup\{\infty\}$ (depending on $\lambda$ and $N$) such that for $m<m_c^{\rm ball}$ balls are the unique minimizers of $E^\gld_{\lambda}(m)$, for $m=m_c^{\rm ball}>0$ balls are minimizers of $E^\gld_{\lambda}(m)$ and for $m>m_c^{\rm ball}$ balls are not minimizers of $E^\gld_{\lambda}(m)$. (This part of \cite{BonCri-14} does not use the assumption $\lambda< N-1$.)

An important result is that $m_c^{\rm ball}>0$, that is, for small $m>0$ balls are minimizers for $E^\gld_{\lambda}(m)$. In the full parameter regime this result is due to \cite{FFMMM-15}, extending earlier results in \cite{KnuMur-13,KnuMur-14,Julin-14,BonCri-14}. The proofs in these papers are based directly or indirectly on the quantitative form of the isoperimetric inequality (see \cite{FuMaPr} and also \cite{FiMaPr,CiLe}) and the regularity theory for quasiminimizers of the perimeter (see, e.g., \cite[Part III]{Ma}). As far as we are aware of, these proofs use compactness arguments and do not give a numerically lower bound on $m_c^{\rm ball}$.

On the other hand, one can show that $m_c^{\rm ball}<\infty$, that is, for large $m>0$ balls are not minimizers for $E^\gld_{\lambda}(m)$. Indeed, setting
$$
m_c^{\rm stab} := \left( \frac{N+1}{\lambda(N-\lambda)}\, \frac{\per B_1}{\frac12\iint_{B_1\times B_1}|x-y|^{-\lambda}\,dx\,dy} \right)^{N/(N-\lambda+1)} |B_1| \,,
$$
one finds that for $m<m_c^{\rm stab}$ the ball is stable against small volume-preserving perturbations and for $m>m_c^{\rm stab}$ it is unstable. (Stability here means that the Hessian is positive definite except for zero modes coming from translations. Instability means that the Hessian is not positive semidefinite.) This computation goes back to Bohr and Wheeler \cite{BoWh} for $N=3$, $\lambda=1$ and can be found in the general case in \cite{BonCri-14,FFMMM-15}. Clearly, $m_c^{\rm ball}\leq m_c^{\rm stab}$, so the former quantity is indeed finite.

\medskip

\emph{Nonexistence.} Let us discuss the nonexistence of minimizers for $E^\gld_{\lambda}(m)$. For fixed $\lambda$ and $N$ we set
$$
m_c^{\rm n.e.} := \sup\left\{ m>0 :\ \text{there is a minimizer for}\ E^\gld_{\lambda}(m) \right\}.
$$
Then, if $\lambda\leq 2$ (and $\lambda<N$, as always), one can show that $m_c^{\rm n.e.}<\infty$, that is, there is no minimizer for large $m$. This is due to \cite{KnuMur-13,KnuMur-14,LuOtt-14,FrNa}. It seems to be unknown whether $m_c^{\rm n.e.}$ is finite or not for $2<\lambda<N$.

In \cite{FraKilNam-16} it is shown that for $\lambda=1$, $N=3$, one has $m_c^{\rm n.e.}\leq 8$. This is to be compared with $m_c^{\rm stab}=10$ for these values of $\lambda$ and $N$. Thus there is a regime $8<m< 10$, where balls are stable local minimizers, but not global minimizers. For comparison, for these values of $\lambda$ and $N$ one has $m_*= 5(2^{1/3}-1)/(1-2^{-2/3})\approx 3.512$.

\begin{problem}
	For $N=3$ and $\lambda=1$, show that balls are minimizers for $m\leq m_*$ and there are no minimizers for $m> m_*$. In which parameter region of $\lambda$'s and $N$'s is the analogous conjecture valid?
\end{problem}

The following two problems are special cases of the previous one.

\begin{problem}
	Do there exist minimizers for $E^\gld_{\lambda}(m)$ for arbitrarily large $m$ in case $2<\lambda<N$?
\end{problem}

\begin{problem}
	Find an explicit numercial lower bound on $m_c^{\rm ball}$, in particular, in the case $N=3$ and $\lambda=1$. 
\end{problem}

We conclude this section by briefly mentioning two further, related models.

The first one concerns the liquid drop model in the presence of a neutralizing background. This problem is motivated, for instance, by the physics of neutron stars and there are interesting mathematical questions; see, e.g., \cite{KnuMurNov-15}. For simplicity we focus here on the case $\lambda=N-2$ in dimension $N\geq 3$, although there are similar versions in dimensions $N=1,2$ \cite{FraLie-19}. For a (large) parameter $L>0$ one sets $\Lambda_L:=(0,L)^N$ and considers the minimization problem
$$
E_L(\rho)\! :=\! \inf\left\{ \per\Omega + \frac12 \!\iint_{\Lambda_L\times\Lambda_L}\!\!\!\!\! \frac{(\1_\Omega(x)-\rho)\,(\1_\Omega(y)-\rho)}{|x-y|^{N-2}}\,dx\,dy : \Omega\subset\Lambda_L, |\Omega|=\rho\, |\Lambda_L| \right\}.
$$
(Sometimes, the kernel $|x-y|^{-N+2}$ is replaced by a constant multiple of the periodic or Neumann Green's function of the Laplacian and the perimeter is replaced by its periodic version or a relative perimeter, but this does not qualitatively change the results discussed below.)

A major open problem is to prove that (for $N=3$, for simplicity) there are $0<\rho_{c1}<\rho_{c2}<1/2$ such that the following holds for minimizers for $E_L(\rho)$ for large $L>0$ `in the bulk': for $0<\rho<\rho_{c1}$, minimizers are periodic with respect to a three-dimensional lattice, for $\rho_{c1}<\rho<\rho_{c2}$, minimizers are periodic with respect to a two-dimensional lattice and for $\rho_{c2}<\rho\leq 1/2$, minimizers are periodic with respect to a one-dimensional lattice. For $1/2<\rho<1$, the situation reverses, with $1-\rho$ replacing $\rho$. This would correspond to what is known as `nuclear pasta' phases in astrophysics.

A fundamental result by Alberti, Choksi and Otto \cite{AlChOt} gives precise bounds on the energy distribution of minimizers that are indicative of the emergence of a regular (e.g., periodic) structure. More precise results about the structure of minimizers are restricted only to the dilute regime. The case $\rho\sim L^{-3}$ is treated in \cite{ChoPel-10} (see also \cite{CicSpa-13} and references therein), the case $\rho\sim L^{-2}$ in \cite{KnuMurNov-15} and the case $\rho\ll 1$ (independently of $L$) in \cite{EmFrKo}.

\medskip

The second generalization of the generalized liquid drop model concerns the addition of an external potential $V$,
\begin{equation*}
	\inf\left\{ \mathcal E_\lambda^\gld[\Omega] + \int_\Omega V\,dx :\ \Omega\subset\R^N \ \text{measurable},\ |\Omega|= m \right\}.
\end{equation*}
Lu and Otto \cite{LuOtt-15} suggested this model with $V(x) = -Z|x|^{-1}$ in $N=3$, $\lambda=1$ as a toy problem for the ionization conjecture in Thomas--Fermi--Dirac--von Weizs\"acker theory and proved that there is no minimizer for $m\geq Z + C \min\{1,Z^{2/3}\}$. Nonexistence for $m\geq Z + C \min\{1,Z^{1/3}\}$, as well as the ionization conjecture in Thomas--Fermi--Dirac--von Weizs\"acker theory was proved in \cite{FraNamvdB}. For more on the ionization conjecture, also for more complicated models, we refer to \cite{Nam}.

\medskip

Finally, returning to the standard liquid drop model with $\lambda=1$ and $N=3$, we mention the open problem to make the global bifurcation picture of Bohr and Wheeler \cite{BoWh} rigorous. For an initial local bifurcation result, see \cite{Fr}.


\section{A simple model for flocking}

In this section we consider the energy functional \eqref{eq:f} and the corresponding minimization problem \eqref{eq:fmin}. We assume throughout that $0<\lambda<N$ and $\alpha>0$.

It is easy to see that there is a minimizer of $E^\f_{\lambda,\alpha}(m)$ for any $m>0$ \cite{CFT}. We would like to understand properties of minimizers and, in particular, qualitative changes in these properties as $m$ varies. For instance, one is interested in the existence of the following three `phases' \cite{FraLie-18}. A first, `liquid' phase occurs when any minimizer $\rho$ for $E^\f_{\lambda,\alpha}(m)$ satisfies $\rho<1$ almost everywhere. A second, `intermediate' phase occurs when there is a minimizer $\rho$ for $E^\f_{\lambda,\alpha}(m)$ such that $\{ 0<\rho<1\}$ has positive measure strictly less than $m$. A third, `solid' phase occurs when any minimizer $\rho$ for $E^\f_{\lambda,\alpha}(m)$ satisfies $\rho=1$ almost everywhere.

\medskip

\emph{Some initial results.}
The case $N\geq 3$, $\lambda=N-2$ and $\alpha=2$ can be solved explicitly \cite{BCT} and one finds that there is an explicit $m_N\in (0,\infty)$ such that the unique (up to translations) minimizer for $E^\f_{\lambda,\alpha}(m)$ is a multiple of the characteristic function of a ball of measure $m_N$ if $m\leq m_N$ and the characteristic function of a ball of measure $m$ if $m>m_N$. In particular, in this special case, the second, intermediary phase does not occur.

In the case $2\leq\alpha\leq 4$ (and any $N\geq 1$ and $0<\lambda<N$), one can show that for any $m>0$ minimizers of $E^\f_{\lambda,\alpha}(m)$ are unique up to translations \cite{Lo} and, in particular, radially symmetric. This relies on an interesting convexity argument. Moreover, the case $N=3$, $\lambda=1$ and $\alpha=4$ is explicitly solved in \cite{Lo}. In particular, there are critical constants $0<m'<m''<\infty$ such that the system is in phase one for $m\leq m'$, in phase two for $m'<m<m''$ and in phase three for $m\geq m''$.

\medskip

\emph{Small $m$ regime.} In \cite{FraLie-18} it is shown that for $N=3$ and $\lambda=1$ (and any $\alpha\geq 1$) there is an $m_*>0$, depending on $\alpha$, such that for all $m<m_*$ any minimizer $\rho$ of $E^\f_{1,\alpha}(m)$ satisfies $\rho<1$ almost everywhere. This result extends, with the same proof, to the case $\lambda=N-2$ in arbitrary dimension $N\geq 3$.

The proof relies on the fact, due to \cite{CDM}, that for $\lambda=N-2$ minimizing measures of the problem
\begin{equation}
	\label{eq:elambdaalpha}
	E_{\lambda,\alpha} := \inf\left\{ \frac12 \iint_{\R^N\times\R^N} \left( |x-y|^{-\lambda} + |x-y|^\alpha\right)d\mu(x)\,d\mu(y) :\ \mu\in P(\R^N) \right\}
\end{equation}
are absolutely continuous with respect to Lebesgue measure with a bounded density. Here $P(\R^N)$ denotes the set of Borel probability measures on $\R^N$. More precisely, one needs a bound on the density depending only on $N$ and $\alpha$.

There are also results in \cite{CDM} concerning the problem $E_{\lambda,\alpha}$ for $0\leq N-2<\lambda<N$ and certain assumptions on $\alpha$. Using these results one should be able to prove that for certain $N$, $\lambda$, $\alpha$, there is an $m_*'>0$, depending on $N$, $\lambda$, $\alpha$, such that for all $m<m_*'$ there are minimizers $\rho$ of $E^\f_{\lambda,\alpha}(m)$ satisfying $\rho<1$.

\medskip

\emph{Large $m$ regime.} Under the assumption $\lambda<N-1$ it is shown in \cite{FrLi} that there is an $m^*<\infty$, depending on $N$, $\lambda$, $\alpha$, such that for $m>m^*$ the only minimizers of $E^\f_{\lambda,\alpha}(m)$ are characteristic functions of balls. The assumption on $\lambda$ is optimal in the sense that for $N-1\leq \lambda<N$ and any $m>0$, balls are not even critical points for the problem $E^\f_{\lambda,\alpha}(m)$.

The results in \cite{FrLi} improve earlier results in \cite{BCT} for $\alpha=2$ and in \cite{FraLie-18} for $\lambda=N-2$, obtained by different methods.

The technique used in \cite{FrLi} is that of symmetric decreasing rearrangement and, more precisely, a quantitative version of the Riesz rearrangement inequality. This quantitative version is due to M.~Christ \cite{Ch}, with some minor extensions and a partially alternate proof in \cite{FrLi2}. As an aside, we mention that from the quantitative Riesz rearrangement inequality one can derive quantitative rearrangement inequalities for Riesz potentials. Those were proved, simultaneously and independently, in a restricted range in \cite{FuPr}; see also \cite{BuCh0,YaYa,BuCh}.

\medskip

Let us conclude this section by mentioning some open problems. Relatively little seems to be known about minimizers of $E_{\lambda,\alpha}^\f(m)$ outside of the asymptotic regimes $m\to 0$ and $m\to\infty$.

\begin{problem}
	Study qualitative properties of minimizers of $E_{\lambda,\alpha}^\f(m)$.
\end{problem}

Concrete questions to be studied are, for instance, the following. Known examples of minimizers are radially symmetric. Can symmetry breaking occur? For arguments in favor of this, see \cite{BCH}. Is the support of a minimizer convex? As $m$ increases, do the regions $\{\rho>0\}$ and $\{\rho=1\}$ increase (fixing the center of mass, for instance), where $\rho$ is a minimizer? Are minimizers concave or convex on their supports for $\alpha<2$ and $\alpha>2$, respectively?

In view of the above small $m$ results, it would be interesting to better understand the case $0<\lambda<N-2$. We consider the minimization problem \eqref{eq:elambdaalpha} and wonder whether the result from \cite{CDM} extends to $0<\lambda<N-2$. An affirmative answer would be related to the existence, for small $m$, of minimizers $\rho$ for $E^\f_{\lambda,\alpha}(m)$ with $\rho<1$ almost everywhere.

\begin{problem}
	For $0<\lambda<N-2$, are minimizers $\mu$ of $E_{\lambda,\alpha}$ absolutely continuous with respect to Lebesgue measure with a bounded density?
\end{problem}

In view of the large $m$ results for $\lambda<N-1$, it seems interesting to investigate in more detail the case $N-1\leq\lambda<N$. We expect that minimizers for large $m$ have values close to one in a large core region and then drop down to zero in a relatively small region. It would be interesting to find the scaling behavior of these regions and, if possible, the transition profile.

\begin{problem}
	For $N-1\leq\lambda<N$ study the shape of minimizers of $E^\f_{\lambda,\alpha}(m)$ for large~$m$.
\end{problem}

\emph{The dynamical problem.} The energy function $\mathcal E^\f_{\lambda,\alpha}$ considered on functions $0\leq\rho\leq 1$ leads via a formal Wasserstein-2 gradient flow to an evolution equation called the constrained aggregation equation; see \cite{CrKiYa,CrTo}. It would be interested to understand the long time behavior of solutions to this equation. In particular, for $\lambda<N-1$ and large $m$ such that characteristic functions of balls are the only optimizers for $E^\f_{\lambda,\alpha}(m)$, one might wonder whether the solution approaches the characteristic function of a ball for large times.


\section{The generalized Keller--Segel model}

In this section we consider the energy functional \eqref{eq:gks} and the corresponding minimization problem \eqref{eq:gksmin}. We assume throughout that $0<q<1$ and $\alpha>0$. We summarize the results from \cite{CDDFH} and \cite{CDFL}.

\medskip

The basic fact is that $E_{q,\alpha}^\gKS=-\infty$ for $0<q\leq N/(N+\alpha)$ and $E_{q,\alpha}^\gKS>-\infty$ for $N/(N+\alpha)<q<1$ \cite[Proposition 20]{CDDFH}. Thus, in the following discussion we will always assume that $q>N/(N+\alpha)$.

It is known and elementary that the case $\alpha=2$ (and any $N/(N+2)<q<1$) can be solved explicitly by expanding the square $|x-y|^2$ and setting the center of mass to zero; see \cite[Corollary 6 and Proposition 20]{CDDFH}. We comment below on the case $\alpha=4$, which can also be solved to some extent.

It is deeper that the case $q=2N/(2N+\alpha)$ can be solved explicitly as well. This was observed by Dou and Zhu \cite{DoZh}, who discovered a conformal symmetry in this case, similarly as in Lieb's work on the Hardy--Littlewood--Sobolev inequality \cite{Li}. The case $q=2N/(2N+\alpha)$ is also of some conceptual importance. If we reinstate the mass in the variational problem \eqref{eq:gksmin} and define $E_{q,\alpha}^\gKS(m)$ in the natural way, then
$$
E_{q,\alpha}^\gKS(m)=m^\frac{2N-(2N+\alpha)q}{N-\alpha - Nq} E_{q,\alpha}^\gKS \,.
$$
Thus, for $q=2N/(2N+\alpha)$, $E_{q,\alpha}^\gKS(m)$ is independent of $m$. As we will see, there are differences between the cases $q>2N/(2N+\alpha)$ and $q<2N/(2N+\alpha)$.

\medskip

\emph{Existence in the superconformal case.} In the case $2N/(2N+\alpha)<q<1$, there is a minimizer for $E_{q,\alpha}^\gKS$ \cite[Proposition 8]{CDDFH}, and any minimizer is radially symmetric with respect to some point, nonincreasing with respect to the distance from this point and positive almost everywhere \cite[Lemma 9]{CDDFH}. Symmetric decreasing rearrangment plays an important role in the proof of existence and in the derivation of the properties of minimizers.

\medskip

\emph{Existence and nonexistence in the subconformal case.} The case $N/(N+\alpha)<q<2N/(2N+\alpha)$ is less understood and there are some open questions about the existence of minimizers. A brief summary of the results in this case is as follows. Either there is a minimizer or there is no minimizer, but instead a generalized minimizer. The latter consists of a symmetric nonincreasing function together with a Dirac delta measure at the center of symmetry. Moreover, sufficient condition for the existence of a `proper' minimizer were given in \cite{CDDFH} and the fact that in some cases there are no minimizers, but only generalized minimizers, was shown in \cite{CDFL}. The existence of a generalized minimizer can be understood as a partial mass concentration phenomenon. We find the appearance of this phenomenon in such a model rather surprising.

Let us be more specific. For $N/(N+\alpha)<q<2N/(2N+\alpha)$, we consider the \emph{relaxed functional}, defined on pairs $(\rho,M)$, where $0\leq\rho\in L^q(\R^N)$ and $M>0$,
\begin{equation}
	\label{eq:rgks}
	\mathcal E_{q,\alpha}^\rgKS[\rho,M] := - \int_{\R^N} \rho^q\,dx + \frac12 \iint_{\R^N\times\R^N} \rho(x) |x-y|^\alpha \rho(y)\,dx\,dy + M \int_{\R^N} |x|^\alpha \rho(x)\,dx \,.
\end{equation}
The corresponding minimization problem is
\begin{equation}
	\label{eq:rgksmin}
	E_{q,\alpha}^\rgKS := \inf\left\{ \mathcal E_{q,\alpha}^\rgKS[\rho]:\ 0 \leq \rho\in L^q(\R^N),\  M>0,\ \int_{\R^N} \rho\,dx +M = 1 \right\}.
\end{equation}
Intuitively, the energy $\mathcal E_{q,\alpha}^\rgKS[\rho,M]$ corresponds to the energy functional $\mathcal E_{q,\alpha}^\gKS$ evaluated at $\rho$ plus a Dirac delta measure of mass $M$ at the origin. Making this intuition rigorous, one finds that \cite[Eq.~(5)]{CDDFH}
$$
E_{q,\alpha}^\rgKS = E_{q,\alpha}^\gKS
$$
and that $E_{q,\alpha}^\gKS$ has a minimizer if and only if $E_{q,\alpha}^\rgKS$ has a minimizer $(\rho_*,M_*)$ with $M_*=0$. Moreover, the same arguments as those applied for $q>2N/(2N+\alpha)$ imply that $E_{q,\alpha}^\rgKS$ has a minimizer \cite[Proposition 10]{CDDFH} and that for any minimizer $(\rho_*,M_*)$ the function $\rho_*$ is radially symmetric with respect to some point, nonincreasing with respect to the distance from this point and positive almost everywhere \cite[Lemma 9]{CDDFH}.

In view of the above discussion, for $N/(N+\alpha)<q<2N/(2N+\alpha)$, the problem of existence of minimizers for $E_{q,\alpha}^\gKS$ is equivalent to the existence of a minimizer $(\rho_*,M_*)$ for the problem $E_{q,\alpha}^\rgKS$ with $M_*=0$. In \cite{CDDFH} we gave sufficient conditions for this. Namely, for $N=1,2$ there is always a minimizer for $E_{q,\alpha}^\gKS$. The same is true for $N\geq 3$ and $\alpha\leq 2N/(N-2)$. If $N\geq 3$ and $\alpha>2N/(N-2)$, this is true provided $q\geq 1-2/N$ \cite[Proposition 11]{CDDFH}.

In \cite{CDFL}, the case $\alpha=4$ was analyzed and an example of a minimizer for $E_{q,\alpha}^\rgKS$ with $M_*>0$ was given. More precisely, it was shown that for $N\geq 6$, the problem $E_{q,4}^\rgKS$ has a minimizer with $M_*>0$ if $q<(N-2)(3N+4)/ ((N+2)(3N))$. Moreover, this result is optimal, in the sense that for $N\geq 6$ and $q\geq (N-2)(3N+4)/ ((N+2)(3N))$, as well as for $N\leq 5$, every minimizer of the problem $E_{q,4}^\rgKS$ has $M_*=0$. The proof is based on a semiexplicit solution.

The paper \cite{CDFL} contains also numerical experiments that are consistent with the appearance of minimizers with $M_*>0$ for $E_{q,4}^\rgKS$. This concentration phenomenon seems to be more pronounced for larger $N$, smaller $q$ and larger $\alpha$.

\begin{problem}
	Prove the existence of a `large' region of parameters $q, \alpha$ for which $E_{q,\alpha}^\rgKS$ has a minimizer $(\rho_*,M_*)$ with $M_*>0$.
\end{problem}

\emph{Uniqueness.} Uniqueness (up to translations) of minimizers, including minimizers of the relaxed functional, is known in two regimes, namely for $2\leq\alpha\leq 4$ and for $\alpha\geq 1$ and $q\geq 1-1/N$ \cite[Theorem 27]{CDDFH}. The first result follows by an small generalization of a proof by Lopes \cite{Lo}, and the latter by the standard tool of displacement convexity in optimal mass transport.

\medskip

\emph{The dynamical problem.} The energy functional $\mathcal E_{q,\alpha}^\gKS$ or, more precisely, its rescaled version
\begin{equation}
	\label{eq:freeenergy}
	-\frac{1}{1-q} \int_{\R^N} \rho^q\,dx + \frac{1}{2\alpha} \iint_{\R^N\times\R^N} \rho(x) |x-y|^\alpha \rho(y)\,dx\,dy
\end{equation}
appears in connection with the aggregation-diffusion equations
\begin{equation}
	\label{eq:aggdiff}
	\partial_t \rho = \Delta \rho^q + \nabla\cdot (\rho \nabla (W * \rho) ) \,,
	\qquad
	W(x) = \alpha^{-1} |x|^\alpha \,.
\end{equation}
Indeed, this time-dependent equation is the formal gradient flow with respect to the Wasserstein-2 distance of the free energy functional \eqref{eq:freeenergy}. Minimizers or, more generally, critical points of the free energy functional, restricted to probability densities, should play an important role for the long time behavior of solutions of \eqref{eq:aggdiff}. It seems particularly interesting to investigate whether in the dynamical setting there is a concentration effect similar to what we have seen for minimizing sequences for $E_{q,\alpha}^\gKS$ in case there is no minimizer, or equivalently there is a minimizer for $E_{q,\alpha}^\rgKS$ with $M_*>0$.

\begin{problem}
	Investigate the long time behavior of solutions of \eqref{eq:aggdiff} in the case where $E_{q,\alpha}^\rgKS$ has a minimizer with $M_*>0$.
\end{problem}

\medskip

To conclude this section, we mention that while we have focused on the free energy functional \eqref{eq:freeenergy} in the case $\alpha>0$ and $0<q<1$, it has been studied for all $q>0$ and $\alpha>-N$. (Here we use the convention that $\alpha^{-1} |x-y|^\alpha$ is understood as $\ln |x-y|$ for $\alpha=0$ and $(1-q)^{-1} \rho^q$ is understood as $- \rho\ln\rho$ for $q=1$.) The nonexistence phenomenon via partial mass concentration that we discussed above, however, appears at most in the region $\alpha>0$ and $0<q<1$. The case $\alpha>0$ and $q\geq 1$ is treated in \cite[Appendix B]{CDDFH}. For $N=2$ and $q=1$ and $\alpha=0$ one obtains the original Keller--Segel free energy functional.


\appendix

\section{The generalized liquid drop model in 1D}\label{sec:app}

In this appendix we consider the minimization problem $E^\gld_\lambda(m)$ in the generalized liquid drop model for $0<\lambda<1$ in dimension $N=1$. We will show that for $m\leq m_*$, single intervals are the unique (up to sets of measure zero) minimizers and for $m>m_*$ there are no minimizers. The computations are elementary.

It is well-known (see, e.g.,~\cite[Proposition 12.13]{Ma}) that any set in $\R$ of finite measure and finite perimeter coincides, up to sets of measure zero, with a finite number of bounded intervals with disjoint closures. Moreover, the perimeter is twice the number of intervals. Clearly, if there are more than one interval, these intervals want to be infinitely far apart. Therefore,
\begin{align*}
	E^\gld_\lambda(m) & = \inf\left\{ 2K + \frac12 \sum_{k=1}^K \int_{-m_k/2}^{m_k/2} \int_{-m_k/2}^{m_k/2} \frac{dx\,dy}{|x-y|^\lambda}:\ K\in\N\,,\ \sum_{k=1}^K m_k = m \right\} \\
	& = \inf\left\{ 2K + \frac{1}{(1-\lambda)(2-\lambda)} \sum_{k=1}^K m_k^{2-\lambda} : \ K\in\N\,,\ \sum_{k=1}^K m_k = m \right\} \\
	& = \inf_{K\in\N} \left( 2K + \frac{1}{(1-\lambda)(2-\lambda)} K^{-1+\lambda} m^{2-\lambda} \right)
\end{align*}
and there is a minimizer if and only if the infimum occurs at $K=1$. Here we used
$$
\sum_{k=1}^K m_k^{2-\lambda} \geq K^{-1+\lambda} \left( \sum_{k=1}^K m_k \right)^{2-\lambda}
$$
(with equality if and only if all $m_k$ are equal). The infimum is attained at $K=1$ if and only if $2+ (1-\lambda)^{-1}(2-\lambda)^{-1} m^{2-\lambda}\leq 2K + (1-\lambda)^{-1}(2-\lambda)^{-1} K^{-1+\lambda} m^{2-\lambda}$ for all $K\geq 2$, which is the same as
$$
m \leq \left( 2(1-\lambda)(2-\lambda) \inf_{K\geq 2} \frac{K-1}{1-K^{-1+\lambda}} \right)^{1/(2-\lambda)} = \left( \frac{2(1-\lambda)(2-\lambda)}{1-2^{-1+\lambda}} \right)^{1/(2-\lambda)} = m_* \,.
$$
Here we used the fact that $\kappa\mapsto (\kappa-1)/(1-\kappa^{-1+\lambda})$ is increasing on $(1,\infty)$. This proves the claimed result.


\bibliographystyle{amsalpha}

\end{document}